\documentclass[10pt]{amsart}

\usepackage{amsfonts}

\usepackage{graphicx}
\usepackage{amsmath}
\usepackage{setspace}

\newtheorem{theorem}{Theorem}
\newtheorem{assumption}{Assumption}

\newtheorem{condition}[theorem]{Condition}
\newtheorem{conjecture}[theorem]{Conjecture}
\newtheorem{corollary}[theorem]{Corollary}

\newtheorem{definition}[theorem]{Definition}
\newtheorem{example}[theorem]{Example}

\newtheorem{lemma}[theorem]{Lemma}

\newtheorem{proposition}[theorem]{Proposition}
\newtheorem{remark}[theorem]{Remark}

\usepackage[english]{babel}
\usepackage[left=3cm,right=3cm,top=2.5cm,bottom=3cm]{geometry}

\usepackage{amsmath,amssymb,amsthm,bbm,color,graphics,version}
\usepackage{mathrsfs}

\newcommand{\bearno}{\begin{eqnarray*}}
\newcommand{\enarno}{\end{eqnarray*}}

\newtheorem{open}{Open Problem}[section]
\newtheorem{conj}{Conjecture}[section]

\title{Tail behaviour of the area under a random process, with applications to queueing systems, insurance and percolations}

\author{Rafa{\l} Kulik}
\thanks{email: rkulik@uottawa.ca; University of Ottawa, Department of Mathematics and Statistics, 585
King Edward Av., K1N 6N5 Ottawa, ON, Canada}
\author{Zbigniew Palmowski}
\thanks{email: zpalma@math.uni.wroc.pl (corresponding author);
Wroc{\l}aw University, Mathematical Institute, Pl. Grunwaldzki 2/4,
50-384 Wroc{\l}aw, Poland}

\date{\today}

\begin{document}

\begin{abstract}
The areas under workload process and under queuing process in a single server queue over the busy period
have many applications not only in queuing theory but also in risk theory or percolation theory.
We focus here on the tail behaviour of distribution of these two integrals.
We present various
open problems and conjectures, which are supported by partial results for some
special cases.
\end{abstract}

\maketitle

\pagestyle{myheadings} \markboth{\sc R.\ Kulik --- Z.\ Palmowski}
{\sc Tail behaviour of the area under a random process}

\vspace{1.8cm}

\tableofcontents

\newpage

\section{Introduction}
In the past two decades there has been done enormous amount of work
on asymptotics for queueing systems. Tail behaviour of steady state
queueing process $\{Q(t),t\ge 0\}$, workload $\{W(t),t\ge 0\}$ or
busy period $\tau$ in standard systems has been well-understood in
both light and heavy tailed case. Surprisingly, however, very little
is known on tail behaviour of integral functionals of the form
\begin{equation}\label{eq:model}
I_f(T):=\int_0^{T} f(X(u))\, du,
\end{equation}
where $\{X(t),t\ge 0\}$ is a stochastic process (typically, $X=Q$ or
$X=W$), $f$ is a deterministic function and $T$ is either $\tau$ or
deterministic (finite or infinite). Such integrals appear naturally
in analysis of ATM. The reader is referred to references given in
\cite{BBP} and \cite{KP}. Recently, in \cite{Antunes}, the authors
connected mean bit rate in time varying $M/M/1$ queue with moments
of the integral $\int_0^{\tau}Q(u)\, du$ in a corresponding standard
$M/M/1$ system.

However, applications of such integrals go beyond queueing systems.
Let $\{S(t),t\ge 0\}$ be a standard risk process. Integrals
$\int_0^T1_{\{S(u)<0\}}S(u)\, du$, where $T$ is deterministic (i.e.
integrated negative part of the risk process), are suggested in
\cite{Loisel} as possible risk measures. Further extensions are
given in a multivariate setting. Furthermore, as in
\cite{Kearney2004}, integrals $\int_0^{\tau}Q(u)\, du$ in
$Geo/Geo/1$ queue and (as a limit) in $M/M/1$ system have particular
interpretation in compact percolations. Last but not least, if $X$
is L\'{e}vy process, integrals $\int_0^{\infty}\exp(-X(u))\, du$
have applications in financial mathematics, see \cite{MaulikZwart}.
Another applications are coming from the actuarial science, where
very often regulated processes are considered and integral
functionals from a regulation random mechanism are investigated.

\section{Subexponential asymptotics}\label{sec:subexp}
Consider a stable $GI/GI/1$ queue. Denote by $\{T,T_i,i\ge 0\}$ and
$\{S,S_i,i\ge 0\}$ two stationary i.i.d. and mutually independent
sequences of interarrival and service times, respectively. Let
$\lambda_T=1/{\rm E}[T]$, $\lambda_S=1/{\rm E}[S]$ and
$\rho=\lambda_T/\lambda_S$. Let $\{Q(t),t\ge 0\}$ be a stationary
queueing process and $\tau=\inf\{t\geq 0: Q(t)=0\}$ the corresponding busy period. We shall
assume that the distribution $F$ of service time $S$ is
subexponential (denotes as $F\in {\mathcal S}$).
The distribution $G$ is subexponential
when
$$
  \lim_{x\to\infty}\overline{G^{*2}}(x)/\overline{G}(x) = 2,
$$
where $G^{*2}$ is the convolution of $G$ with itself and
$\overline{G}$ denotes the tail distribution given by
$\overline{G}(x)=1-G(x)$.

Heuristically, the large area $\int_0^{\tau}Q(u)\, du$ is realized
by a customer with large service time $S_0$, say, who arrives at the
very beginning of the busy period and blocks the server. During that
time $S_0$, according to the Law of Large Numbers (LLN),
approximately $\lambda_TS_0$ customers arrive and the queue length
process increases linearly. Hence the area under queueing process
before reaching maximum is asymptotically equivalent to the area of
a triangle: $\frac{1}{2\lambda_T}\bar{Q}(\tau)^2$, where
$\bar{Q}(\tau)$ is maximum of the queue length process over the busy
period (see Figure 1). After passing the maximum, the queue length
process behaves according to the LLN, decreasing almost linearly to
$0$ with the slope $\lambda_S-\lambda_T$. Thus, the area under the
queueing process after reaching the maximum is equivalent to
$\frac{1}{2(\lambda_S-\lambda_T)}\bar{Q}(\tau)^2$ and $\bar Q(\tau)$
is equivalent to $\rho(\lambda_T-\lambda_S)\tau$. This heuristic
leads to the following conjecture.

\begin{figure}
\setlength{\unitlength}{0.5mm}
\begin{picture}(150,150)

\put(30,55){$\lambda_T$} \put(170,55){$\lambda_T-\lambda_S$}
\put(90,130){$\bar{Q}(\tau)$}

\put(10,15){\line(0,0){10}} \put(25,18){\line(0,0){4}}
\put(193,18){\line(0,0){4}}

\put(10,5){0}
\put(193,5){$\tau$} 
\put(10,20){\vector(4,0){200}}
\put(10,20){\vector(0,0){120}}\put(15,120){$Q(t)$}

\put(10,23){\line(4,0){15}} \put(25,23){\line(1,2){50}}
\put(75,123){\line(1,-1){8}} \put(83,115){\line(1,2){5}}
\put(88,125){\line(1,-1){105}}

\end{picture}\caption{A typical behaviour of a heavy tailed queueing process.}
\end{figure}
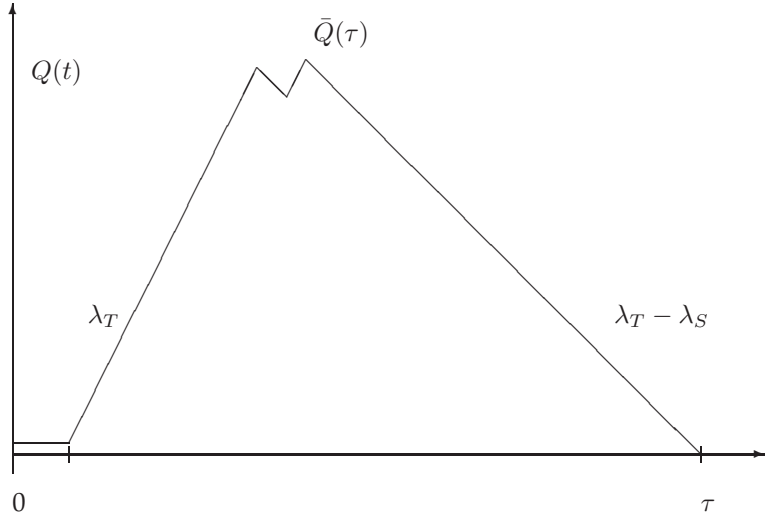
\begin{conj}
If $F \in {\mathcal S}$, then
\begin{equation}\label{eq:2}
P\left(\int_0^{\tau}Q(u)\, du>x\right)\sim
P\left(\tau> \sqrt{\frac{2x}{\rho(\lambda_S-\lambda_T)}}\right).
\end{equation}
\end{conj}
Similarly, for the workload process the heuristic is as follows. The
most likely way for the area to be large is that one early big
service time occurs and apart from this, everything in the cycle
develops normally. Using LLN and ignoring random fluctuations, this
leads to the conclusion that the workload goes to zero with negative
rate $-(1-\rho)$. Thus the area exceeds level $x$ iff the area of
the triangle with the sides $\tau(1-\rho)$ and $\tau$ is greater
than $x$, hence when
$$\frac{1}{2}\tau^2(1-\rho)>x\ ,$$
which suggests the following conjecture.
\begin{conj}
If $\bar F \in {\mathcal S}$, then
\begin{equation}\label{eq:1}
P\left(\int_0^{\tau}W(u)\, du>x\right)\sim
P\left(\tau>\sqrt{\frac{2x}{1-\rho}}\right).
\end{equation}
\end{conj}
The statements (\ref{eq:2}) and (\ref{eq:1}) were proven in
\cite{KP} and \cite{BBP}, respectively, under regularly varying
assumption of the service time, that is $\bar F(x)=x^{-\alpha}L(x)$,
where $\alpha>1$ and $L$ is slowly varying at infinity. Furthermore,
in \cite{KP} one needs additionally that
\begin{equation}\label{interlight}
\lim_{t\to\infty}\frac{t^{1+\varsigma}P(T>t)}{\bar F(t)}=0
\end{equation}
holds with $\varsigma>0$.

The tail behaviour of the busy period $\tau$ can be identified in
terms of $\bar F$ for a large subclass of ${\mathcal S}$:
\begin{equation}\label{eq:busy-period}
P(\tau>x)\sim \frac{1}{1-\rho}\bar F((1-\rho)x),
\end{equation}
see \cite{BaltrunasDaleyKlup}, \cite{JelMom} and \cite{Zwartbusy}. In
particular, using (\ref{eq:busy-period}) in the regularly varying
case, together with (\ref{eq:2}) and (\ref{eq:1}) yields exact
asymptotics for area under queueing process and workload,
respectively (see \cite{KP} and \cite{BBP}).

\section{Light-tailed asymptotics}

As in Section \ref{sec:subexp} we consider a stable $GI/GI/1$ queue.
Here, we assume that the service time $S$ is light-tailed, that is
there exists $\theta>0$ such that ${\rm E}[\exp(\theta S)]<\infty$.

Under the above assumptions for the queueing process, we have the
following open problem:
\begin{open}\label{open:1}
Find exact asymptotics of
$$P\left(\int_0^{\tau}Q(u)\, du>x\right).$$
\end{open}

We suppose that
$$P\left(\int_0^{\tau}Q(u)\, du>x\right)\sim
Cx^{-1/4}\exp(-\psi \sqrt{x})
$$
for some constants $\psi$ and $C$.

We
suggest above asymptotics, believing that it should be the same like
for $M/M/1$ queue, which was found in \cite{GuilPin} under two
conjectures on p. 391 and it is in the following form:
\begin{equation}\label{mm1q}
P\left(\int_0^{\tau}Q(u)\, du>x\right)\sim
\frac{1-\rho}{\rho\sqrt{2\pi\psi}}x^{-1/4}\exp(-\psi \sqrt{x}),
\end{equation}
where
$$
\psi=2\sqrt{-2(1-\rho)+(1+\rho)\log \rho}.
$$
In the proof authors used the Laplace transform method. We are
not aware of any probabilistic proof of this result and we do not
know if Conjectures 1 and 2 in \cite{GuilPin} hold true.

Unfortunately, we have not managed to produce any heuristic for this
result either. The idea of the piecewise linear most likely
trajectory seems to produce wrong expression. In particular, define
the new probability measure $\tilde P$:
\begin{equation}\label{eq:5}
\frac{d\tilde P_{|{\mathcal F}_n}}{dP_{|{\mathcal F}_n}}=e^{\gamma
\sum_{i=1}^n(S_i-T_{i-1})},
\end{equation}
where ${\mathcal F}_n=\sigma(T_1,S_1,\ldots,T_{n},S_n)$ and $\gamma$ solves the equation ${\rm E}[\exp (\gamma(T-S))]=1$.
Let $\tilde\rho=\frac{\tilde {\rm
E}[T]}{\tilde {\rm E}[S]}>1$.
Consider the most likely path coming from large deviation theory for large cycle maxima, that is
trajectory that develops along the line with the slope $\tilde{\rho}-1$ and after getting maximum
behaves 'normally', that is goes to zero linearly
with negative rate $-(1-\rho)$.
This trajectory produces wrong asymptotics for $M/M/1$ queue since by Kyprianou \cite{Kyprfather} we have then:
$$P\left(\int_0^{\tau}Q(u)\, du>x\right)\sim P\left(\tau>\sqrt{2\frac{1+\rho}{1-\rho}x}\right)\sim
Cx^{-3/4}e^{-\gamma \sqrt{2\frac{1+\rho}{1-\rho}}\sqrt{x}}$$ for
some constant $C$ and $\gamma=(1-\sqrt{\rho})^2\mu$. The explanation
might come from papers \cite{Meyn0,Meyn1, Meyn2}, where the optimal
trajectories in a sense of large deviation theory for the mean value
of a reflected random walk $W_n$ are considered. The relationship
between area under the queue length and sum of $W_n$ is clear when
observing queue length process at arrival and departure epochs. This
papers suggest that optimal path, though still concave in general,
might not be piecewise linear. In fact, this papers suggest that the
optimal trajectory $q(t)$ of the queueing process for the large
value of the area on the cycle should solve the following equation:
\begin{equation}\label{eq:meyn}
\nabla
I\left(\frac{d}{dt}q(t)\right)=\lambda^*(\tau-t),\end{equation}
where $\lambda^*>0$ and $I$ is a rate function for the increment
process. Still, solving (\ref{eq:meyn}) explicite seems to be a
difficult task and transferring it into finding the asymptotic tail
distribution of $\tau$ and the integral even more cumbersome. This
problem shows also the need of simulations of the area
$\int_0^{\tau}Q(u)\, du$ to check if (\ref{eq:meyn}) indeed produces
non-piecewise-linear optimal path. In general, we are not aware of
any simulations results concerning area under a random process.

Similar open problem one can pose for $P\left(\int_0^{\tau}W(u)\, du>x\right)$.

\section{Extensions}
\subsection{Transient case}
Let us consider the $M/M/1$ queue with $\rho=1$. In
\cite{Kearney2004} the author established the following asymptotics
for the critical case $\rho=1$:
\begin{equation}\label{eq:4}
P\left(\int_0^{\tau}Q(u)\,
du>x\right)\sim\frac{3^{1/3}}{\Gamma(1/3)\psi_0\psi_1}x^{-1/3},
\end{equation}
where $\psi_0$ and $\psi_1$ are explicit constants.
\begin{open}
Consider a transient $GI/GI/1$ queue such that $\bar F(x)$ is
regularly varying. What is the tail asymptotics of
$\int_0^{\tau}Q(u)\, du$ and $\int_0^{\tau}W(u)\, du$ ?
\end{open}
As mentioned above in the Introduction, this type of questions
should have connections with heavy-tailed critical percolations.
\subsection{Utility functions, discounting}
Consider again a stable $GI/GI/1$ queue. Assume that $f$ is a
deterministic function. We are interested in the tail behaviour of $
\int_0^{\tau}f(W(u))\, du $. In an insurance context, $f(\cdot)$ may
play a role of the utility function. Of course, one can formulate
the corresponding problem for the queuing process, where $f$ gives the costs of maintaining the system.\\

The heuristic given in the subexponential case in Section
\ref{sec:subexp} suggests that studying the asymptotics of latter integral is
equivalent to study the behaviour of
$$
P\left(\int_0^{\tau}f((\tau-u)(1-\rho)) \, du  >x \right),
$$
when $x$ becomes large.
For example, if one considers $f(x)=x^k$, $k\ge 0$, it leads to the
following conjecture.
\begin{conj}\label{conj:1}
Assume that $\bar F\in {\mathcal S}$ and $f(x)=x^k$. Then,
$$
P\left(\int_0^{\tau}f(W(u))\, du >x\right)\sim P\left(\tau
>\sqrt[k+1]{\frac{(k+1)x}{(1-\rho)^k}}\right).
$$
\end{conj}
However, this is not clear for us what should be expected for a very
general function $f$.\\

Now, assume that $\theta>0$ and consider the process $X(t)=\exp(-\theta
t)f(W(t))$, $t\ge 0$. Similarly, the corresponding integral can be
interpreted in an insurance context: $\exp(-\theta t)$ is the
discounting factor and $f(\cdot)$ is the utility function.

The subexponential heuristic leads us again to
$$
P\left(\int_0^{\tau}\exp(-\theta u)f((\tau-u)(1-\rho)) \, du  >x
\right).
$$
For further heuristic, consider $f(x)=x^k$. Then
\begin{eqnarray*}
\lefteqn{\int_0^{\tau}\exp(-\theta u)f((\tau-u)(1-\rho)) \,
du=(1-\rho)^k\int_0^{\tau}\exp(-\theta u)(\tau-u)^k\, du}\\
&=& (1-\rho)^k\exp(-\theta\tau)\int_0^{\tau}\exp(\theta s)s^k\, ds
\\
&=& (1-\rho)^k\exp(-\theta\tau)
\left[\frac{1}{\theta}\exp(\theta\tau)\tau^k-\frac{k}{\theta}\int_0^{\tau}\exp(\theta
s)s^{k-1}\, ds\right].
\end{eqnarray*}
Integrating further by parts the last integral we see that the
leading term will be $O(\exp(\theta \tau)\tau^{k-1})$. This leads to
the following conjecture.
\begin{conj}
Assume that $F\in {\mathcal S}$ and $f(x)=x^k$. Then,
$$
P\left(\int_0^{\tau}\exp(-\theta u)f(W(u))\,  du>x\right) \sim
P\left(\tau
>\sqrt[k]{\frac{\theta x}{(1-\rho)^k}}\right).
$$
\end{conj}
This conjecture should be compared with Conjecture \ref{conj:1}, showing how
discounting leads to the change in the asymptotic behaviour.\\

Let us consider now a light tailed case. The large area for the
integral $\int_0^{\tau}W(u)\; du$ is strictly connected to a large
$\bar W(\tau)$. It does not seem to be the case if one considers
$\int_0^{\tau}\exp(-\theta u)f(W(u))\, du$, i.e. large values of
$W(u)$ may be killed by the discount factor.
\begin{open}
Under the conditions of Open Problem \ref{open:1}, find the
asymptotics for
$$
P\left(\int_0^{\tau}\exp(-\theta u)f(W(u))\,  du>x\right).
$$
\end{open}
\subsection{Finite horizon}\label{sec:finite}
All the conjectures and open problems above may be formulated in
case of $\int_0^T$, where $T$ is finite and deterministic. In
particular, if $S(t)=v+ct-\sum_{i=1}^{N(t)}S_i$, $t\ge 0$, is the
classical risk insurance process, then moments of
$$
\int_0^TS(u)1_{\{S(u)<0\}}\, du
$$
are proposed in \cite{Loisel} as particular risk measures. In that
paper the authors established the asymptotics for the expected value
of the integral, when $T$ is fixed and the initial capital $v$
becomes large. Clearly, this problem can be re-formulated for the
dual workload process. Furthermore, all the above mentioned
extensions (discounting, utility functions) may be considered.

\subsection{Multivariate case}
As in Section \ref{sec:finite}, let us consider insurance context.
We assume that a company has two lines of business and customers
arrive according to a renewal process $N(t)$, $t\ge 0$, with generic interarrival time $T$. The $i$th
customer has a claim $(S_{i,1},S_{i,2})$. For example, a car
accident may cause a claim for driving and liability insurance
(see for example \cite{APPins}, \cite{Loisel} and references therein). If
two lines of business are considered separately, we are interested
in the tail behaviour of
$$
\left(\int_{0}^TS_1(u)1_{\{S_1(u)<0\}}\, du ,
\int_{0}^TS_2(u)1_{\{S_2(u)<0\}}\, du\right),
$$
where $S_1$ and $S_2$ are risk process associated with the
corresponding lines of business. One can consider the tail behaviour
when $T$ is fixed and a vector of initial capital $(x_1,x_2)$
becomes large (in a particular sense, usually dependence is linear, that is $x_1=ax_2$ for some fixed $a>0$).

In queuing theory one can consider e.g. parallel queues, where customers are coming into the system according to renewal process
$N(t)$ and each service time (we assume here that all of them are i.i.d.) is proportionally divided into two servers (see e.g. \cite{LishoutMichel1} and \cite{LishoutMichel2}).
In this case, $S_{i,1}=bS_{i,2}$ for some fixed $b$. Let $W_i(t)$ and $Q_i(t)$ will be the workload and queueing process on $i$th server ($i=1,2$).
The busy period $\tau$ is understood here as the minimum of the busy periods on
both servers.
We could analyze the tail behaviour of the following "tail" probability:
$$P\left(\int_0^\tau W_1(u)\,du>x\quad\mbox{and}\quad \int_0^\tau W_2(u)\,du >ax\right)\qquad \mbox{as $x\to\infty$}, $$
under assumption that distribution $F$ of $S_{i,1}$ is subexponential or light-tailed.
Similar considerations could be analyzed for the bivariate occupation process.

\subsection{Regulated processes}

Apart of original process $X(t)$ in (\ref{eq:model}), one can
consider regulated process:
$$U(t)=X(t)-L(t),$$
where $L(t)$ is left-continuous, increasing process, adapted with
respect to the natural filtration of $X$. For example, if $X$ is a
spectrally negative L\'{e}vy risk process, then $L(t)$ might be
cumulative dividends payment up to time $t$ payed according to some
strategy. The most often used strategy is so-called barrier strategy
in which all surpluses above a given level $a$ are transferred
(possibly subject to a discount rate) to a beneficiary. In this case
$L(t)=a\vee \sup_{s\leq t}X(s)-a$ and for given utility function we
are interested in in the following random variable:
\begin{equation}\label{calkarefl1}I_f(\tau)=\int_0^\tau e^{-\theta u}f(L(u))\,du,\end{equation}
(possibly with $\theta=0$), where $\tau=\inf\{t\geq 0: U(t)<0\}$ is a ruin time.

We can also analyze discounted cumulative dividends payed up to ruin
time:
\begin{equation}\label{calkarefl2}\int_0^\tau e^{-\theta u}\,dL(u),\end{equation}
which is equivalent to finding the tail asymptotics of integral:
$$\int_0^\infty e^{-\theta Z_1(t)}\mathbf{1}_{\{\Delta Z_2(u)<a\}}\,du,$$
where $(Z_1(t),Z_2(t))$ is a bivariate subordinator and $\Delta
Z_2(u)=Z_2(u)-Z_2(u-)$ is a size of jump of $Z_2$ at time $u$ (see
\cite{Al, APP, KPal, Sch} for details and other references). In the
case if $a=\infty$ we end up with integral from exponential function
of L\'{e}vy process, see e.g. \cite{Bert1, Bert2, Carmona,
MaulikZwart} for this kind of functional. Note that in particular
$I=\int_0^\infty e^{-\theta Z_1(t)}\,du$ solves the following
equation: $I\stackrel{{\mathcal D}}{=}\int_0^\tau e^{-\theta Z_1(t)}\,du+e^{-\theta
Z_1(\tau)}I$ and analyzing $I$ is strongly related then with
properties of $e^{-\theta Z_1(\tau)}$ (see e.g. \cite{Goldie,
Kesten}).

With this problem there is also related integral with respect to the general L\'{e}vy process:
$$\int_0^TY(u)\,dX(u)$$
for appropriate predictable integrand $Y$ and fixed $T$ (see
\cite{HultLindskog} for multivariate regularly varying setting).\\

In queueing systems with single server and finite capacity $a$, $L(t)$ corresponds then to cumulative lost
of work. If only proportion of information is lost, then we can consider e.g. $L(t)=c\mathbf{1}_{\{U(t)>a\}}$
for some $c>0$. In this case formally so-called refracted process $U$ solves the following stochastic equation:
$$dU(t)=dX(t)-c\mathbf{1}_{\{U(t)>a\}}\,dt;$$
see \cite{Andronnie}. In this case we can also try to find
asymptotic tail of (\ref{calkarefl1}) and (\ref{calkarefl2}).

\section*{Acknowledgements}
The authors thank Sean P. Meyn for inspiring discussions and
suggestion of optimal path (\ref{eq:meyn}). This work is partially
supported by the Ministry of Science and Higher Education of Poland
under the grants N N201 394137 (2009-2011) and N N201 525638
(2010-2011).



\begin{thebibliography}{99}
\bibitem{Al}
H. Albrecher and S. Thonhauser, Optimality Results for Dividend Problems in Insurance, RACSAM Rev. R. Acad. Cien. Serie A. Mat., 103(2) (2009), 295--320.
\bibitem{Antunes}
N. Antunes, C. Fricker, F. Guillemin, and P. Robert, Perturbation analysis of the area swept under the queue length process
of a variable $M/M/1$ queue. Performance 2005.
\bibitem{APPins}
F. Avram, Z. Palmowski, Z. and M. Pistorius,
Exit problem of a two-dimensional risk process from a cone:
exact and asymptotic results, Annals of Applied Probability, 18(6) (2008), 2421--2449.
\bibitem{APP}
F. Avram, Z. Palmowski, Z. and M. Pistorius, On the optimal
dividend problem for a spectrally negative L\'{e}vy process,
Ann. Appl. Probab., 17 (2007),  156--180.
\bibitem{BaltrunasDaleyKlup}
A. Baltrunas, D.J. Daley, and C. Kl\"{u}ppelberg, Tail behaviour of
the busy period of a $GI/GI/1$ queue with subexponential service
times, Stoch. Proc. Appl., 111 (2004)
237--258.
\bibitem{Bert1}
J. Bertoin and M. Yor, On the entire moments of self-similar Markov processes and exponential
functionals of L\'{e}vy processes, Ann. Fac. Sci. Toulouse Math. (6), 11(1) (2002), 33--45.
\bibitem{Bert2}
J. Bertoin, P. Biane, and M. Yor, Poissonian exponential functionals, q-series, q-integrals,
and the moment problem for log-normal distributions, In Seminar on Stochastic Analysis,
Random Fields and Applications IV, volume 58 of Progr. Probab. (2004), 45--56. Birkh\H{a}user,
Basel.
\bibitem{Loisel}
R. Biard, S. Loisel, C. Maccib and N. Veraverbeke, Asymptotic
behavior of the finite-time expected time-integrated negative part
of some risk processes and optimal reserve allocation,
\textit{Preprint.}
\bibitem{BBP}
A.A. Borovkov, O.J. Boxma and Z. Palmowski, On the integral of the
workload process of the single server queue, J. Appl. Probab., 40 (2003), 200--225.
\bibitem{Carmona}
P. Carmona, F. Petit, and M. Yor,  On the distribution and asymptotic results for exponential
functionals of L\'evy processes, In Exponential functionals and principal values related to
Brownian motion, Bibl. Rev. Mat. Iberoamericana (1997), 73--130.
\bibitem{DenisovShneer}
D. Denis Denisov and V. Shneer, Asymptotics for first passage times
of L\'{e}vy processes and random walks. Preprint at arxiv 0712.0728,
(2007).
\bibitem{Meyn0} K. R. Duffy and S. P. Meyn, Most likely paths to error when estimating the mean of a
reflected random walk, Performance Evaluation, 67(12) (2010), 1290--1303.
\bibitem{Goldie}
C. M. Goldie, Implicit renewal theory and tails of solutions of random equations. Ann. Appl.
Probab., 1(1) (1991), 126--166.
\bibitem{GuilPin}
F. Guillemin and D. Pinchon, On the area swept under the occupation
process of an $M/M/1$ queue in a busy period, Queueing Systems, 29
(1998), 383--398.
\bibitem{HultLindskog}
H. Hult and F. Linskog, Extremal behavior of stochastic integrals driven by regularly varying L\'{e}vy processes,
Ann. Probab., 35 (2007), 309--339.
\bibitem{JelMom}
P. R. Jelenkovi{\'c} and P. Mom{\v{c}}ilovi{\'c}, Large deviations
of square root insensitive random sums, Math. Oper. Res., 29 (2004)
398--406.
\bibitem{Kearney2004}
M. J. Kearney, On a random area variable arising in discrete-time
queues and compact directed percolation, Journal of Physics. A.
Mathematical and General, 37 (2004), 8421--8431.
\bibitem{Kesten}
H. Kesten, Random difference equations and renewal theory for products of random matrices,
Acta Math., 131 (1973), 207--248.
\bibitem{Kyprfather}
E.K. Kyprianou,
On the quasi-stationary distribution of the virtual waiting
time in queues with Poisson arrivals, J.~Appl.~Probab., 8 (1971), 494--507.
\bibitem{KPal}
A. Kyprianou and Z. Palmowski,
Distributional study of De Finetti's dividend problem
for a general L\'{e}vy insurance risk process,
J. Appl.Probab., 44(2) (2007), 428--443.
\bibitem{Andronnie}
A. Kyprianou and R. Loeffen, Refracted L\'{e}vy processes,
Ann. Inst. H. Poincar\'{e} Probab. Statist., 46 (2010), 24--44.
\bibitem{KP}
R. Kulik and Z. Palmowski, Tail behaviour of the area under queue
length process of a single server queue with regularly varying
service times, Queueing Systems, 50 (2005), 299--323.
\bibitem{LishoutMichel1}
P. Lieshout and M. Mandjes Asymptotic analysis of Levy-driven tandem queues, Queueing Systems, 60 (2008), 203--226.
\bibitem{LishoutMichel2}
P. Lieshout and M. Mandjes, Brownian tandem queues, Mathematical Methods in Operations Research, 66 (2007), 275--298.
\bibitem{MaulikZwart}
K. Maulik and B. Zwart, Tail asymptotics for exponential functionals
of {L}\'evy processes, Stoch. Proc. Appl., 116 (2006),
156--177.
\bibitem{Meyn1} S. P. Meyn, Large deviation asymptotics and control variates for simulating large functions,
Ann. Appl. Probab., 16(1) (2006), 310--339.
\bibitem{Meyn2} S. P. Meyn, Control Techniques for Complex Networks, Cambridge University Press, Cam-
bridge, 2007.
\bibitem{PalmowskiRolski}
Z. Palmowski and  T. Rolski,
On the exact asymptotics of the busy period in GI/G/1 queues,
Adv. Appl. Probab., 38(3) (2006), 792--803.
\bibitem{Sch}
H. Schmidli, Stochastic Control in Insurance, Springer, New York, 2008.
\bibitem{Zwartbusy}
A.P. Zwart, Tail asymptotics for the busy period in the $GI/G/1$
queue, Math. Oper. Res., 26 (2001), 485--493.
\end{thebibliography}
\end{document}